# A Note on the Existence of the Multivariate Gamma Distribution

Thomas Royen

Fachhochschule Bingen, University of Applied Sciences

e-mail: thomas.royen@t-online.de

**Abstract.** The $p$ - variate gamma distribution in the sense of Krishnamoorthy and Parthasarathy exists for all positive integer degrees of freedom $\nu$ and at least for all real values $\nu > p-2,\ p \geq 2$. For special structures of the "associated" covariance matrix it also exists for all positive $\nu$. In this paper a relation between central and non-central multivariate gamma distributions is shown, which implies the existence of the $p$ - variate gamma distribution at least for all non-integer $\nu$ greater than the integer part of $(p-1)/2$ without any additional assumptions for the associated covariance matrix.

## 1. Introduction

The p-variate chi-square distribution (or more precisely: "Wishart-chi-square distribution") with $\nu$ degrees of freedom and the "associated" covariance matrix $\Sigma$ (the $\chi_p^2(\nu, \Sigma)$ - distribution) is defined as the joint distribution of the diagonal elements of a $W_p(\nu, \Sigma)$ - Wishart matrix. Its probability density (pdf) has the Laplace transform (Lt)

$$|I_p + 2\Sigma T|^{-\nu/2}, \tag{1.1}$$

with the $(p \times p)$ - identity matrix $I_p$, $\nu \in \mathbb{N}$, $T = diag(t_1,...,t_p)$, $t_j \geq 0$, and the associated covariance matrix $\Sigma$, which is assumed to be non-singular throughout this paper.

The $p$ - variate gamma distribution in the sense of Krishnamoorthy and Parthasarathy [4] with the associated covariance matrix $\Sigma$ and the "degree of freedom" $\nu = 2\alpha$ (the $\Gamma_p(\alpha, \Sigma)$ - distribution) can be defined by the Lt

$$|I_p + \Sigma T|^{-\alpha} \tag{1.2}$$

of its pdf

$$g_\alpha(x_1,...,x_p; \alpha, \Sigma). \tag{1.3}$$

For values $\nu = 2\alpha \in \mathbb{N}$ this distribution differs from the $\chi_p^2(\nu, \Sigma)$ - distribution only by a scale factor 2, but in this paper we are only interested in positive non-integer values $2\alpha$ for which $|I_p + \Sigma T|^{-\alpha}$ is the Lt of a pdf and not of a function also assuming any negative values. These values $\alpha$ are called here "admissible values". The admissibility of all values $2\alpha > p-1$ follows from the existence of the $W_p(2\alpha, \Sigma)$ - distribution, and the admissibility of $2\alpha \in (p-2, p-1)$ follows from formula (1.4) below. Smaller values of $\alpha$ are admissible at least with some additional assumptions for $\Sigma$. Sufficient and necessary conditions for $\Sigma$ entailing infinite divisibility

---

of the Lt $|I_p + \Sigma T|^{-1}$ (i.e. all $\alpha > 0$ are admissible) are found in [1] and [3]. According to [1], $|I_p + \Sigma T|^{-1}$ is infinitely divisible if and only if there exists any signature matrix $S = diag(s_1, ..., s_p)$, $s_j = \pm 1$, for which $S\Sigma^{-1}S$ is an M-matrix. For the infinite divisibility of a more general class of multivariate gamma distributions see [2].

At least all values $2\alpha > m-1$, $m \in \mathbb{N}$, $m < p$, are admissible for "$m$-factorial" covariance matrices $\Sigma = W^{-2} + AA^T$ ($\Leftrightarrow W\Sigma W = I_p + BB^T$, $B = WA$ with rows $b^j$) with a suitable matrix $W = diag(w_1, ..., w_p)$, $w_j > 0$, and a real $(p \times m)$-matrix $A$ of the lowest possible rank $m$. This follows from the representation

$$g_\alpha(x_1, ..., x_p; \alpha, \Sigma) = E\left(\prod_{j=1}^{p} w_j^2 g_\alpha(w_j^2 x_j, \tfrac{1}{2} b^j S b^{jT})\right) \tag{1.4}$$

of the $\Gamma_p(\alpha, \Sigma)$-pdf (see [7] and [8]), where

$$g_\alpha(x, y) = e^{-y} \sum_{n=0}^{\infty} g_{\alpha+n}(x) \frac{y^n}{n!}$$

is the non-central gamma density with the non-centrality parameter $y$ and the central gamma densities $g_{\alpha+n}$, and the expectation refers to the $W_m(2\alpha, I_m)$-Wishart matrix $S$. With $W^{-2} = \lambda I_p$, where $\lambda$ is the lowest eigenvalue of $\Sigma$, it follows $m \leq p-1$. A special case – entailing infinite divisibility – is a one-factorial $\Sigma = W^{-2} + aa^T$ with a real column $a$.

In the following section it will be shown that all values $2\alpha > [(p-1)/2]$ (the integer part of $(p-1)/2$) are admissible without any further assumptions for $\Sigma$. This result is obtained as a corollary of a relation between central and non-central multivariate gamma densities given in theorem 1. This relation was already derived in a similar form for integer values $\nu = 2\alpha$ in [6], but the non-integer values require a different proof by means of the Lt.

## 2. A Sufficient Condition for the Existence of the $\Gamma_p(\alpha, \Sigma)$-Distribution

Let be given any partition

$$\Sigma = \begin{pmatrix} \Sigma_{11} & \Sigma_{12} \\ \Sigma_{21} & \Sigma_{22} \end{pmatrix} \tag{2.1}$$

of the non-singular $(p \times p)$-covariance matrix $\Sigma$ with $(p_i \times p_i)$-matrices $\Sigma_{ii}$ and set

$$\Sigma_0 = \Sigma_{11} - \Sigma_{12} \Sigma_{22}^{-1} \Sigma_{21}, \tag{2.2}$$

which is a non-singular covariance matrix too.

For $2\alpha > \max(p_1 - 1, p_2 - 1)$ the Wishart $W_{p_2}(2\alpha, \Sigma_{22})$-pdf and the non-central $W_{p_1}(2\alpha, \Sigma_0, \Delta)$-pdf exist where the symmetrical positive semi-definite $(p_1 \times p_1)$-matrix $\Delta$ is any "non-centrality matrix" of rank $k \leq p_1$. The diagonal of a random matrix $Z$ has a non-central $\Gamma_{p_1}(\alpha, \Sigma_0, \Delta)$-distribution if $2Z$ has a $W_{p_1}(2\alpha, \Sigma_0, \Delta)$-distribution, which exists (apart from integer values $2\alpha$ with $\text{rank}(\Delta) \leq \min(p_1, 2\alpha)$) for all $\alpha > (p_1 - 1)/2$ with $\text{rank}(\Delta) \leq p_1$ (see [5]).

The corresponding $\Gamma_{p_1}(\alpha, \Sigma_0, \Delta)$ - pdf

$$g(x_1, ..., x_{p_1}; \alpha, \Sigma_0, \Delta) \tag{2.3}$$

has the Lt

$$|I_1 + \Sigma_0 T_1|^{-\alpha} \operatorname{etr}(-T_1(I_1 + \Sigma_0 T_1)^{-1}\Delta) \tag{2.4}$$

with $T_1 = diag(t_1, ..., t_{p_1})$. (In the literature the "non-centrality matrix" is frequently defined in a different non-symmetrical way.) Let $\mathcal{C}_{p_2}$ denote the set of all $(p_2 \times p_2)$ - correlation matrices $C$ and let

$$2Y = 2X^{1/2} C X^{1/2} \tag{2.5}$$

be a $W_{p_2}(2\alpha, \Sigma_{22})$ - matrix with a random $C \in \mathcal{C}_{p_2}$ and $X = diag(X_{p_1+1}, ..., X_p)$, where the elements $X_j$ have the gamma-densities $\sigma_{jj}^{-1} g_\alpha(\sigma_{jj}^{-1} x_j) = \sigma_{jj}^{-\alpha} x_j^{\alpha-1} \exp(-\sigma_{jj}^{-1} x_j)/\Gamma(\alpha)$ with the $\sigma_{jj}$ from the diagonal of $\Sigma_{22}$. The density of $Y$ is given by

$$(\Gamma_{p_2}(\alpha))^{-1} |\Sigma_{22}|^{-\alpha} |Y|^{\alpha-(p_2+1)/2} \operatorname{etr}(-\Sigma_{22}^{-1} Y) \tag{2.6}$$

with the multivariate gamma function $\Gamma_{p_2}(\alpha) = \pi^{p_2(p_2-1)/4} \prod_{j=1}^{p_2} \Gamma\left(\alpha - \frac{j-1}{2}\right)$.

Now, with the notations from (2,1), (2.2), (2.3), (2.5), (2.6), we show for the $\Gamma_p(\alpha, \Sigma)$ - pdf from (1.3):

**Theorem 1.** If $p = p_1 + p_2$ and $2\alpha > \max(p_1 - 1, p_2 - 1)$, then

$$g_\alpha(x_1, ..., x_p; \alpha, \Sigma) = (\Gamma_{p_2}(\alpha))^{-1} |\Sigma_{22}|^{-\alpha} \times$$

$$\int_{\mathcal{C}_{p_2}} g_\alpha(x_1, ..., x_{p_1}; \alpha, \Sigma_0, \Sigma_{12}\Sigma_{22}^{-1} X^{1/2} C X^{1/2} \Sigma_{22}^{-1} \Sigma_{21}) |X|^{\alpha-1} |C|^{\alpha-(p_2+1)/2} \operatorname{etr}(-\Sigma_{22}^{-1} X^{1/2} C X^{1/2}) dC.$$

**Remark.** The simple special case with $p_2 = 1$, $C = 1$ (and $\Sigma_{22} = 1$) was already given by theorem 2 in [11].

As a corollary we obtain

**Theorem 2.** The function $g_\alpha(x_1, ..., x_p; \alpha, \Sigma)$ with the Lt $|I_p + \Sigma T|^{-\alpha}$ from (1.2) is a $\Gamma_p(\alpha, \Sigma)$ - pdf at least for $2\alpha \in \mathbb{N} \cup ([(p-1)/2], \infty)$.

**Proof of theorem 2.** Choose for $p_1$ or $p_2$ in theorem 1 the value $[(p+1)/2]$. □

**Proof of theorem 1**. The equation in theorem 1 will be verified by the Lt of both sides. The left side has the Lt $|I_p + \Sigma T|^{-\alpha}$ and we get with the Schur complement for determinants

$$|I_p + \Sigma T| = \begin{vmatrix} I_1 + \Sigma_0 T_1 + \Sigma_{12}\Sigma_{22}^{-1}\Sigma_{21} T_1 & \Sigma_{12} T_2 \\ \Sigma_{21} T_1 & I_2 + \Sigma_{22} T_2 \end{vmatrix} =$$

$$|I_2 + \Sigma_{22}T_2| \cdot |I_1 + \Sigma_0 T_1 + \Sigma_{12}\Sigma_{22}^{-1}\Sigma_{21}T_1 - \Sigma_{12}T_2(I_2 + \Sigma_{22}T_2)^{-1}\Sigma_{21}T_1| =$$

$$|I_2 + \Sigma_{22}T_2| \cdot |I_1 + \Sigma_0 T_1 + \Sigma_{12}\Sigma_{22}^{-1}(I_2 - \Sigma_{22}T_2(I_2 + \Sigma_{22}T_2)^{-1})\Sigma_{21}T_1| =$$

$$|I_2 + \Sigma_{22}T_2| \cdot |I_1 + \Sigma_0 T_1 + \Sigma_{12}\Sigma_{22}^{-1}(I_2 + \Sigma_{22}T_2)^{-1}\Sigma_{21}T_1| =$$

$$|I_1 + \Sigma_0 T_1| \cdot |I_2 + \Sigma_{22}T_2| \cdot |I_1 + \Sigma_{12}\Sigma_{22}^{-1}(I_2 + \Sigma_{22}T_2)^{-1}\Sigma_{21}T_1(I_1 + \Sigma_0 T_1)^{-1}|. \tag{2.7}$$

With (2.4), (2.5) and (2.6) we find for the right side the Lt

$$|I_1 + \Sigma_0 T_1|^{-\alpha} (\Gamma_{p_2}(\alpha))^{-1} |\Sigma_{22}|^{-\alpha} \times$$

$$\int_{Y>0} \text{etr}(-T_1(I_1+\Sigma_0 T_1)^{-1}\Sigma_{12}\Sigma_{22}^{-1}Y\Sigma_{22}^{-1}\Sigma_{21})|Y|^{\alpha-(p_2+1)/2}\,\text{etr}(-\Sigma_{22}^{-1}Y - T_2 Y)\,dY =$$

$$|I_1 + \Sigma_0 T_1|^{-\alpha} (\Gamma_{p_2}(\alpha))^{-1} |\Sigma_{22}|^{-\alpha} \times$$

$$\int_{Y>0} |Y|^{\alpha-(p_2+1)/2}\,\text{etr}\left(-(\Sigma_{22}^{-1}\Sigma_{21}T_1(I_1+\Sigma_0 T_1)^{-1}\Sigma_{12}\Sigma_{22}^{-1} + \Sigma_{22}^{-1} + T_2)Y\right) dY =$$

$$|I_1 + \Sigma_0 T_1|^{-\alpha} |\Sigma_{21}T_1(I_1+\Sigma_0 T_1)^{-1}\Sigma_{12}\Sigma_{22}^{-1} + I_2 + \Sigma_{22}T_2|^{-\alpha},$$

(see e.g. formula (2.2.6) in [12]), and

$$|I_1 + \Sigma_0 T_1| \cdot |I_2 + \Sigma_{22}T_2| \cdot |I_2 + \Sigma_{21}T_1(I_1+\Sigma_0 T_1)^{-1}\Sigma_{12}\Sigma_{22}^{-1}(I_2+\Sigma_{22}T_2)^{-1}| =$$

$$|I_1 + \Sigma_0 T_1| \cdot |I_2 + \Sigma_{22}T_2| \cdot |I_1 + \Sigma_{12}\Sigma_{22}^{-1}(I_2+\Sigma_{22}T_2)^{-1}\Sigma_{21}T_1(I_1+\Sigma_0 T_1)^{-1}|,$$

which coincides with formula (2.7), where the last identity follows from the general equation

$$\begin{vmatrix} I_1 & A_{12} \\ -B_{21} & I_2 \end{vmatrix} = |I_2 + B_{21}A_{12}| = |I_1 + A_{12}B_{21}|. \quad \Box$$

**Some further remarks.** If $\Sigma = W^{-2} + AA^T$ is $m$-factorial, where $A$ – and consequently $B = WA$ – may contain any mixture of real or pure imaginary columns, then the function with the representation from (1.4) has again the Lt $|I_p + \Sigma T|^{-\alpha}$ and it is a $\Gamma_p(\alpha, \Sigma)$-pdf at least for all values $2\alpha > m-1 \geq [(p-1)/2]$. For special structures of $\Sigma$ smaller values of $2\alpha$ are possible, e.g. with an $m$-factorial $\Sigma = W^{-2} + AA^T$ with a real matrix $A$ and $m-1 < [(p-1)/2]$. Furthermore, let be, e.g., $p_2 < p_1$, $\Sigma_0 = W_0^{-2} + A_0 A_0^T$ with a real $(p_1 \times m_0)$-matrix $A_0$ of rank $m_0$ and $m_{12} = \text{rank}(\Sigma_{12}) \leq p_2$. Then, at least $2\alpha > \max(m_0 + m_{12} - 1, p_2 - 1)$ is admissible and $\max(m_0 + m_{12} - 1, p_2 - 1) < [(p-1)/2]$ is possible for low values of $m_0$ and $m_{12}$.

On the other hand it is at present an open question if there exist some $(p \times p)$-covariance matrices $\Sigma$ for which $2\alpha$ is inadmissible for some values $2\alpha \in ([(p-3)/2], [(p-1)/2])$, $p \geq 5$.

A consequence of theorem 2 is the extension of the inequality

$$G_p(x_1,...,x_p;\alpha,\Sigma) > G_{p_1}(x_1,...,x_{p_1};\alpha,\Sigma_{11})G_{p-p_1}(x_{p_1+1},...,x_p;\alpha,\Sigma_{22}),$$

$$x_1,...,x_p > 0, \text{ rank}(\Sigma_{12}) > 0, \tag{2.8}$$

for the $p$ - variate cumulative $\Gamma_p(\alpha,\Sigma)$ - distribution function $G_p$. This inequality was proved for $2\alpha \in \mathbb{N}$ and for all values $2\alpha > p-2$ in [9] (see also [10]), and it implies the famous Gaussian correlation inequality for $2\alpha = 1$. Now, this inequality can be extended to all non-integer values $2\alpha > [(p-1)/2]$ without any further assumptions for $\Sigma$.